\documentclass[a4paper,10pt]{amsart}

\usepackage{algorithm,algorithmic,color,amsmath}

\newtheorem{theorem}{Theorem}[section]

\newtheorem{lemma}[theorem]{Lemma}

\newtheorem{proposition}[theorem]{Proposition}
\newtheorem{corollary}[theorem]{Corollary}

\theoremstyle{definition}
\newtheorem{definition}[theorem]{Definition}

\newtheorem{remark}[theorem]{Remark}

\newtheorem{example}[theorem]{Example}

\numberwithin {equation}{section}

\newcommand {\D}  {{\mathcal D}}

\newcommand {\LL}  {{\mathcal L}}

\newcommand {\Vc}  {{\mathcal V}}

\newcommand {\NN}  {{\mathbb N}}

\newcommand {\ZZ}  {{\mathbb Z}}

\newcommand{\kommentar}[1]{} 
\newcommand{\comment}[1]{}

\newcommand{\abs}[1]{\left\vert#1\right\vert}
\newcommand{\set}[1]{\left\{#1\right\}}

\DeclareMathOperator{\Card}{ Card }

\def\vectz #1 #2 {\left(
\begin{array}{c}
 #1 \\ #2
\end{array}\right)}
\def\vektor #1 #2 {\left(
\begin{array}{c}
 #1 \\ \vdots \\ \vdots \\ #2
\end{array}\right)}

\def\msg{{\rm msg }}

\title[Digital semigroups]{Digital semigroups} 

\begin{document}

\author{Horst~Brunotte}
\address{Haus-Endt-Stra{\ss}e 88,  D-40593 D\"usseldorf, Germany\\e-mail: brunoth@web.de}

\date{\today}

\begin{abstract}
The well-known expansion of rational integers in an arbitrary integer base different
from $0, 1, -1$ is exploited to study relations between numerical monoids and certain subsemigroups of the multiplicative semigroup of nonzero integers.
\end{abstract}



\maketitle

{\rm 
\begin{small}
Article in RAIRO -- Theoretical Informatics and  Applications  50  (2016) 67 -- 79.\\
The original publication is available at www.edpsciences.org/ita.\\
DOI: https://doi.org/10.1051/ita/2016005
\end{small}}

\vspace{15mm}

\section*{Introduction}

Recently, Rosales,  Branco and Torr{\~a}o \cite{rosbrantorr} investigated
  sets of positive integers and their relations to the number of decimal
 digits. More precisely, they introduced and thoroughly studied digital semigroups
 which are defined as follows.
 A digital semigroup $D$ is a subsemigroup of the semigroup $(\NN \setminus \set{0}, \cdot)$
such that for all $d\in D$ the set $ \set{n\in \NN \;:\; \ell  (n) = \ell (d) }$ is contained in  $D$; here $\NN$ is the set of nonnegative rational integers and    $\ell  (n)$
denotes the number of digits of $n$ in the usual decimal expansion.
Among other things, the smallest digital semigroup containing a set of positive
integers is determined, and for this purpose a bijective map $\theta$ between the set of digital semigroups
and a certain subset $\LL$ of numerical monoids, namely LD-semigroups, is constructed.
 Recall that a  numerical monoid is a submonoid  of $(\NN, +)$ whose complement in $\NN$
 is finite, and  an LD-semigroup $S$ is a numerical monoid such that there exists a digital
semigroup $D$ with the property   $S=\set{\ell  (d) \;:\; d \in D } \cup  \set{0}$.
It is shown that $\LL$ is a Frobenius variety and that the elements of $\LL$ can
be arranged in a tree. Moreover, LD-semigroups  are characterized by the fact that the
minimum element in each interval of nongaps belongs to the minimal set of generators.
Finally, it is  observed that certain combinatorial configurations introduced
by Bras-Amor\'os and Stokes \cite{brasamorosstokes} are in fact LD-semigroups.

It is well-known that every positive integer can be represented in an arbitrary
integer  base larger than one. Expansions of integers in negative integer bases have apparently been introduced
by Gr\"unwald~\cite{Gr} and rediscovered by several authors; the reader is referred to 
Knuth~\cite{knuthsa}  for more
details. In view of these facts   we extend the notions of digital semigroups and
LD-semigroups coined by Rosales,  Branco and Torr{\~a}o for
decimal expansions to expansions of integers in an arbitrary integer base, i.e.,
instead of the base  $b=10$ we consider an integer base $b \ne  0,1, -1$.
Consequently, we replace the digit set   $ \set{0,1,\ldots, 9}$ by the
canonically chosen set $ \set{0,1,\ldots, \abs{b}-1}$ and simply
apply the prefix $b$ (subscript $b$, respectively) at appropriate places; clearly, by
omitting  $b$  the original notions are recovered.

It turns out that for positive base $b$ essentially  all results coincide with the respective
results presented by  Rosales,  Branco and Torr{\~a}o; however, for negative base $b$
some  modifications have to be taken. In particular, 
 bijective maps $\theta_b$ between the set of certain $b$-digital semigroups
and  specified subsets of $\LL$ play an important role here.

\section{$b$-digital semigroups}
 
In this article we always let  $b \in \ZZ \setminus \set{-1,0,1}$ and denote by $N_b:=  \set{0,1,\ldots, \abs{b}-1}$ the set of all nonnegative integers less than $\abs{b}$.
It  is well-known  that for positive $b$ every positive  integer $z$ can uniquely be represented in the form
\begin{equation}\label{repab}
 z= \sum_{i=0}^{n} u_i b^i \qquad  \qquad (u_0,\ldots, u_n \in  N_b,\; u_n \ne 0);
\end{equation}
similarly, if $b$ is negative then every non-zero integer $z$ can uniquely be written in the form~\eqref{repab}.
Putting\footnote{Obviously, this and some other notions in the sequel depend only on the
sign of $b$. However, our notion facilitates subsequent formulations.} $Z_b:=\NN \setminus \set{0}$ for $b>0$ \;  ($Z_b:=\ZZ \setminus \set{0}$ for $b<0$,
respectively) the  positive integer  
$$\ell_b (z):= n +1$$
is called the length of the representation of $z\in Z_b$ in base $b$, and we consistently set 
$\ell_b (0):=1$. 
Thus,  for every $z\in  Z_b  \cup \set{0}$ the  integer 
$\ell_b (z)$ denotes the number of digits  of the representation of $z$ in base $b$.
Some elementary properties of the length function are collected in  the last section.

\bigskip

We now generalize the fundamental notion of a digital semigroup in the sense explained in the introduction. Further,  for these new objects we  present some examples and properties which will be used in the sequel. 

\begin{definition}
A $b$-digital semigroup $D$ is a subsemigroup  of $(Z_b, \cdot)$
such that $ \Delta_b ( \ell_b (d)) \subseteq D$ for all $d\in D$.
Here we introduce the notation
 $$\Delta_b (n):= \set{z \in Z_b  \;:\;  \ell_b (z) = n }\qquad (n\in \NN\setminus \set{0}).$$
 \end{definition}

\bigskip

Following  \cite{rosbrantorr} we let
$$L_b (A):= \set{\ell_b (a) \;:\; a \in A }$$
for the set $A \subseteq Z_b$, and we apply the commonly used abbreviation
$$\set{z_1, \ldots, z_k, \to }:=\set{z_1, \ldots, z_k }\cup \set{z \in  \ZZ\;:\; z >z_k }$$
for integers $z_1<  \cdots<  z_k$.

\bigskip

Before listing some properties of $b$-digital semigroups we present several examples.
In particular, these examples show that the analogue of \cite[Proposition~2]{rosbrantorr} does not
hold unrestrictedly.

\begin{example}\label{exbigdsgr}  
\begin{enumerate} 
 \item  Let  $D: = \set{1}$ be the trivial subgroup of $(\ZZ\setminus \set{0}, \cdot)$. 
 If $\abs{b}= 2$ then $D$ is a $b$-digital semigroup; however, 
 $L_b (D)$ is not additively closed.  
Trivially,  if $\abs{b}> 2$
then $D$ is not a $b$-digital semigroup.
\item   The set   $ Z_b  \setminus  N_b$
 is a $b$-digital semigroup,   and 
$L_b ( Z_b  \setminus  N_b)= \set{2, \to}$ is a subsemigroup of $(\NN, +)$.
\item   Let  $b< -1, \; \ell_0 \ge 3$ and $D :=  \set{d \in \ZZ\setminus \set{0} \;:\; \ell_b (d) \text{ odd}, \; \ell_b (d) \ge \ell_0}.$
 Then $D\subset \NN\setminus \set{0}$ by Proposition \ref{modlen} below,  $D$ is a $b$-digital semigroup by Lemma~\ref{lacn}, but 
$$L_b (D)= \set{2n+1 \;:\; n \in \NN, \; n \ge (\ell_0  - 1)/ 2 }$$ is not additively closed.
 \item The set
 $$D :=  \set{z  \in \ZZ \;:\; \ell_{-2} (z) \ge 3} =  \ZZ \setminus (\Delta_{-2} (1)\cup \Delta_{-2} (2) \cup  \set{0})$$
is a $(-2)$-digital semigroup, 
and  $L_{-2} (D)=\{3, \to\}$ is additively closed.
\end{enumerate}
 \end{example}
  
\bigskip

The essential ideas for  the proof of the following statements are taken from \cite[Proposition~2]{rosbrantorr}.

\begin{lemma} \label{dsnsgv}
 Let    $D$ be a $b$-digital semigroup.
 \begin{enumerate} 
 \item   If   $x \in L_b (D)$ and $u\in N_b\setminus \set{0}$ then  $u b^{x-1} \in D$.
  \item  If  $x, y \in L_b (D)$  then $x+ y-1 \in L_b (D) \,.$
\item There  exist $x,y \in L_b (D)$ such that
 $\gcd (x, y)=1\,.$
 \item  Let   $x,y\in L_b (D) $.
  If  $b \ge 3$ then $x+y \in L_b (D) $, and if  $b \le -3$ then $x+y+1 \in L_b (D) $.
  \end{enumerate} 
     \end{lemma}
\begin{proof}
(i) By definition we have  $\ell_b (u b^{x-1})=x\in L_b (D)$, hence $u b^{x-1} \in D$.\\
(ii)  By  (i) we have  $ b^{x-1}, \,  b^{y-1}\in D$, hence 
$  b^{x+y-2}\in D$ which yields
$$x+y-1 =\ell_b (b^{x+y-2}) \in L_b (D).$$
(iii)  Pick $x\in L_b (D)$ such that $x>0$. By (ii) we have $y:=2x-1 \in L_b (D)$, and clearly  $\gcd (x, y)=1\,.$\\
(iv)  Pick  $u, v \in N_b$ such that $\abs{b} \le uv < 2 \abs{b}$. Then there exists 
$w \in N_b$ such that 
$$uv =\abs{b} +w\,.$$
By (i) we have  $u b^{x-1}, \, v b^{y-1}\in D$, hence 
$$d:= (\abs{b} +w)b^{x+y-2}=uv  b^{x+y-2}=(u b^{x-1})(v b^{y-1})\in D\,.$$
If $b>0$  we deduce 
$$x+y=\ell_b ( b^{x+y-1})=\ell_b (b\cdot  b^{x+y-2})=\ell_b (d)\in L_b(D),$$
and if $b<0$  we have 
$$x+y+1=\ell_b ( b^{x+y})=\ell_b (b^2\cdot  b^{x+y-2})=\ell_b (((\abs{b}-1) b + b^2)\cdot  b^{x+y-2})=\ell_b (d)\in L_b(D),$$
since
$$\abs{b} +w= b^{2}+ (\abs{b}-1)b + w\,.$$
\end{proof}
 
\bigskip 

Our interest concerns the structure of the set of the lengths of the $b$-adic representations of 
the elements of a $b$-digital semigroup.

\begin{proposition}\label{dsns} 
Let  $D$ be a $b$-digital semigroup.
Then  $L_b (D)\cup \set{0}$ is a numerical monoid provided that  one of the following conditions holds.
\begin{enumerate}
\item  $L_b (D)$ is additively closed.
  \item  $b \ge 3$.
\item $b = 2\; $ and  $ \; 2 \, \cdot \, \min (L_2 (D))  \in L_2 (D)$.
\item  For all $n, m\in \NN$
 the relation $b^n, \,  b^m \in D$ implies $b^{n +m+1} \in D$.
\end{enumerate}
\end{proposition}
\begin{proof}
Set $S:=L_b (D)\cup \set{0}$.\\
(i) Pick $x \in S  \setminus \set{0}$. Then  
Lemma~\ref{dsnsgv}  yields  $2 x-1 \in S$. In view of $\gcd (x, 2x-1)=1$
our assertion now follows  from \cite[Lemma~1]{rosbrantorr}.\\
 (ii) Lemma  \ref{dsnsgv} shows that $S$ is additively closed,
and   then (i)  implies our assertion. \\
(iii) Let $n,m\in S\setminus \set{0}$. 

\medskip

Case 1 \hspace{3cm}  $n =  1$ \quad  or  \quad  $m =  1$

\medskip

Then we have $  \min (S\setminus \set{0}) = 1\in D$. By assumption this yields
$\ell_2 (d) =2$ for some $d\in D$, thus $2\in D$ and further $2^k\in D$
for all $k\in \NN$. But then we have $S=\NN$, and we are done.

\medskip

Case 2 \hspace{3cm}  $n, m  >  1$

\medskip

In view of  Proposition~\ref{merhrp} we have 
$2^{n}-1, \,  2^{m}-1 \in D$,  thus
 $$ d:= (2^{n}-1)  \, (2^{m}-1) \in D. $$
We easily check
 $$ 2^{n+m-1} \le d < 2^{n+m} \,,$$
 and we conclude
 $$n+m =\ell_2 (d) \in S$$
by   Proposition~\ref{merhrp}, and again we are done by (i).\\
(iv) Clear by Lemma~\ref{dsnsgv} and (i). 
 \end{proof}

\section{$b$-LD-semigroups}

In this section  we adapt the notion of an LD-semigroup introduced in  \cite{rosbrantorr}.
We characterize $b$-LD-semigroups and construct a correspondence between
$b$-digital semigroups and $b$-LD-semigroups. Further, several examples and properties of
$b$-LD-semigroups for negative $b$ are listed.
 
\begin{definition}
Let $S$ be a submonoid of $(\NN, +)$. We call  $S$  a $b$-LD-semigroup
 if there exists a
 $b$-digital semigroup $D$ such that  $S= L_b (D) \cup  \set{0}$.
 \end{definition}
 
Now we are in a position to extend  \cite[Theorem~4]{rosbrantorr}
and provide the crucial  characterization of  $b$-LD-semigroups.
For ease of notation, we put $ E_b:=\set{-1}$ for $b> 1 $ and  $E_b:=  \set{-3, -1, 1}$ for $b<-1 $.
 
\begin{theorem}\label{ldsgrchar}
 Let  $S$ be a submonoid of $(\NN, +)$. Then the following statements are
 equivalent:
  \begin{enumerate} 
 \item   $S$ is  a $b$-LD-semigroup.
   \item  $S\ne  \set{0}$ and   $s+t+ e\in S$ for all
  $s,  t\in S\setminus \set{0, 1}$ and $e \in E_b \,.$
  \end{enumerate} 
\end{theorem} 
\begin{proof}
(i) $\implies$ (ii) 
Let  $D$ be a $b$-digital semigroup such that  $S= L_b (D) \cup  \set{0}$. Then we clearly have
 $S\ne  \set{0}$.  Let $s, t  \in S\setminus \set{0, 1}$ and $e \in E_b \,.$ By 
 Lemma~\ref{lacn} there exist $a, c \in Z_b$ such that $s= \ell_b (a), \, t= \ell_b (c)$
 and $\ell_b (a c) = s + t +  e \,.$
 By the properties of  $D$ we know that  $a, c \in D$, thus
$$ s+t+ e = \ell_b (a c)  \in  L_b (D) \subset S\,.$$
(ii) $\implies$ (i)
Since  $S\ne  \set{0}$ the set
$$D_b:=  \set{z\in Z_b \;:\; \ell_b (z) \in S}$$
is nonempty, and we immediately convince ourselves that $S=L_b (D_b) \cup \set{0}.$ 
By construction we have $ \Delta_b ( \ell_b (d)) \subseteq D_b$ for all $d\in D_b$.
Therefore we are left to show
that  $D_b$ is multiplicatively closed.

Let   $a, c \in D_b$, thus $s: = \ell_b (a), \, t: = \ell_b (c) \in S$.
If  $s=1 \text{ or } t=1$ then  $\NN \subseteq S$, and we are done.
 Therefore we may assume  $s, t  >  1$.
If $b> 1$ then our prerequisites and 
Lemma~\ref{lacn} yield some  $e \in  \set{-1, 0}$
such that
\begin{equation}\label{ldsgrchareq}
  \ell_b (a c ) = s+t + e  \in  S.
 \end{equation}
Similarly, if $b < -1$ then there is  some  $e \in  E_b$
such that \eqref{ldsgrchareq} holds.
Thus, in both cases we have shown  $a c \in D_b$.
 \end{proof}

\bigskip

Let us list some direct consequences of this result.

\begin{corollary}\label{ldsgrcharc}
 Let  
 $S$ be a $b$-LD-semigroup.
\begin{enumerate}
  \item If $b  > 1$ then   $S$ is a $c$-LD-semigroup for all $c > 1$.
\item If $b  < -1$ then   $S$ is a $c$-LD-semigroup for all  $ c \in \ZZ \setminus \set{-1,0,1}$.
\end{enumerate}
\end{corollary}

\bigskip

\begin{corollary}\label{ldsgrcharclds}
Every  $b$-LD-semigroup is a numerical monoid.
\end{corollary}
\begin{proof} Using  Theorem~\ref{ldsgrchar} the proof is analogous to 
 \cite[Proposition~2]{rosbrantorr} and left to the reader. 
\end{proof}
 
\bigskip

\begin{remark}
\begin{enumerate}
  \item  Let   $b  > 1$ and  $S$ be a $b$-LD-semigroup. Then $S$ need not be a $c$-LD-semigroup for $c < -1$, e.g., consider $S=\set{0,4,7,\to}$.
  \item  Let $ b<-1$,   $D$ be a $b$-digital semigroup,   $n, m \in L_b (D)$.
 Then there do not exist $d, e  \in D$ such that 
$\ell_b (d)=n, \; \ell_b (e)=m $ and  $\ell_b (d e)=n +m$.
Indeed, if $n+m$ is even then either both $n, m$ are odd or both $n, m$ are even.
In any case the product $d e$ is positive, hence  $\ell_b (d e)$ is odd  
(cf.  Proposition~\ref{modlen}). We similarly argue
in the case  $n+m$ odd.
  \end{enumerate}
\end{remark}

\bigskip
 
In view of  Theorem~\ref{ldsgrchar}  we let
$$\LL:= \set{S \text{ submonoid of } \NN\;:\; S\ne  \set{0}, \; s+t-1\in S \text{ for all }  s, t \in S \setminus  \set{0, 1} }$$
 be the set of all  $b$-LD-semigroups for $b > 1$,
and 
$$\LL_- := \{S \text{ submonoid of } \NN\;:\; S\ne  \set{0}, \; s+t-3, \, s+t-1, \, s+t +1\in S$$
$$ \text{ for all }  s, t \in S \setminus  \set{0, 1} \}$$
 be the set of all $b$-LD-semigroups for $b<  -1$. By what we have seen above, 
 $\LL$ coincides which the respective set in \cite[Section~2]{rosbrantorr}. Moreover,
 $\LL_-$ is a proper subset of $\LL$ (see Example~\ref{exllmin} below), and
 by \cite[Proposition~12]{rosbrantorr} the set 
$\LL$ is a Frobenius variety 
which has  been investigated in detail in \cite{rosbrantorr}. Recall that a Frobenius
variety 
is a nonempty set $\mathcal V$ of numerical semigroups with the following properties:
\begin{enumerate}
\item If $S, T \in \mathcal V$, then $S \cap T \in \mathcal V$. 
\item If $S \in \mathcal V$ and $S \ne  \NN$, then $S \cup \set{F(S)} \in \mathcal V$. 
\end{enumerate}
Here,   for $A\subseteq \NN$ such that 
$ \Card \, (\NN \setminus A) < \infty$
we let $F(A)$ denote the Frobenius number of $A$, i.e.,  the greatest integer which does not belong to $A$.

\bigskip

 In view of our remark above,  
  we now mainly concentrate on the subset $\LL_-$ of the Frobenius variety $\LL$.

\bigskip

Some examples which also illustrate  subsequent results seem appropriate.
As usual,  we denote by  $\msg (S)$ the (unique) minimal set of generators of the numerical monoid $S$. 

\begin{example} \label{exllmin}
\begin{enumerate}
  \item Let  $n\in \NN \setminus \set{0}$.  The LD-semigroups   $S_n:=  \set{0,n,\to}$ appear as the left-most
  branch in the tree of LD-semigroups presented in  \cite[Figure~1]{rosbrantorr}; note that $S_n \in \LL_-$ if and only if $n\ne 2$ since $2+2-3\notin S_2$.  Clearly, 
$\msg (S_n)= \set{n, \ldots, 2 n -1}$, and for $n\ge 2$ we have $F(S_n)= n-1$ and 
$S_n  \setminus \set{n}= S_{n+1} \in \LL_-$, but 
$S_n  \setminus \set{2 n - 1}=  \set{0,n, \ldots, 2 n - 2, 2 n,\to}\notin \LL$.
  \item   $<3, 5,7>,   <4, 5,7> \,  \in \LL_-$, but  $  <4, 6,7,9>  \,  \in  \LL \setminus \LL_-$ since 
  $ 4+4-3=5 \notin <4, 6,7,9>$.  
  \item Trivially, we have  $\NN   \in \LL_-$. By (i) we have  $S:=S_3  \in \LL_-$, and we easily check $\msg (S)= \set{3,4,5}$
  and  $S \setminus  \set{5}= \set{0,3, 4, 6,\to}   \notin \LL$,
  since  $3+3-1 \notin S \setminus \set{5}$.  
Further, we have 
$S \cup   \set{F(S)}= S_2  \notin \LL_-$.
We remark in passing that $\LL_-$ is not a Frobenius pseudo-variety (see \cite{robperros} for details).
  \end{enumerate}
 \end{example}

\bigskip

Motivated by the last example  we establish  the following observation.

\begin{proposition} \label{frobenius}
$\LL_- \setminus \set{\set{0, 3, \to}}$ is a Frobenius variety.
 \end{proposition}
\begin{proof}
The proof follows the same lines as \cite[Proposition~12]{rosbrantorr}.
Set  $S_3:=  \set{0, 3, \to}$.
Clearly, $\Vc:= \LL_-  \setminus \set{S_3} \ne  \emptyset $ since $\NN \in \Vc$. 

It is immediate that $S, T \in \Vc$ implies $S \cap  T \in \ \Vc$.
Indeed, $S \cap  T  \in \LL_-$ by  Theorem~\ref{ldsgrchar}, and by \cite[Section~3]{rosbrantorr} the assumption $S \cap  T  = S_3$ implies $S=S_3$ or  $T=S_3$ which is impossible.  

Now, let  $S \in \Vc$ such that   $S \ne \NN$. Note that $2\notin S$, since
otherwise $1=2+2-3\in S$ which we excluded. Let  $e \in  \set{-3,-1, 1} $ and
  $s, t   \in S \cup   \set{F(S)}$ such that   $s, t  > 1$.
  If $s, t   \in S$ then certainly $s+t + e \in S$. Therefore it remains to consider
  the case $F(S) \in  \set{s, t}$. If $F(S) = s$ then  $s>2$ because otherwise $F(S)= 2$ and 
  $S=S_3$ which is impossible. Thus we may assume $s, t  \ge 3$, hence 
   $s+t+e \ge F(S),$    and we are done.
\end{proof}

\bigskip

Applying the ideas of  \cite[Proposition~14]{rosbrantorr} we can  derive the
following result without difficulty.

\begin{proposition} \label{lgmmsg}
Let $S \in \LL_-$ such that  $ 3 \notin S$, and let  $s \in \msg (S)$.
Then $S \setminus \set{s} \in \LL_- $ if and only if 
  $s -1, \; s +1, \; s +3 \in (\NN\setminus S) \cup  \msg (S)$.
 \end{proposition}
\begin{proof}
Note  that $S\cap \set{-1,1, 3} =  \emptyset $ by our prerequisites.

\medskip

Let  $S \setminus \set{s} \in \LL_- $ and   assume
  $s+e \notin (\NN \setminus S) \cup   \msg (S)$
  for some   $e \in  \set{-1,1, 3} $. Then $s +e \in  S \setminus  \msg (S)$
  and there exist $t, r  \in  S$ such that $s+e =t+r$. In view of
  $$t + r -e= s  \notin S \setminus \set{s} $$
   we infer $S \setminus \set{s} \notin \LL_- $ from Theorem~\ref{ldsgrchar}: Contradiction. 

\medskip

Conversely, let  $t, r  \in  S\setminus \set{0, s}$, thus in particular $t, r  \ne 1$.
Using Theorem~\ref{ldsgrchar} again we see
$ t+r - e   \in S $ for each   $e \in  \set{-1,1, 3} $.
The assumption $t+r - e= s$ leads to $s+ e= t+r \in S\setminus \set{s}$
which implies the contradiction
$s+ e \notin (\NN \setminus S) \cup   \msg (S).$
Thus we have shown $t+r \in S\setminus \set{s}$, and we are done by Theorem~\ref{ldsgrchar}.
\end{proof}
 
\bigskip

Analogously as  \cite[Corollary~15]{rosbrantorr} we can formulate:

\begin{corollary} \label{lgmmsgc}
Let $S \in \LL_-$ such that  $ 3 \notin S$, and let  $s \in \msg (S)$
with $s > F(S)$.
Then $S \setminus \set{s} \in \LL_- $ if and only if 
$s-1 \in (\NN\setminus S) \cup  \msg (S)$ and $\; s+1, \; s+3 \in  \msg (S)$.
   \end{corollary}

\begin{remark}
Note that we cannot renounce the assumption  $ 3 \notin S$ in our two last results. Indeed,  choose  $s= 3$ and consider the semigroups
 $ <3, 5,7>  $ for  Proposition~\ref{lgmmsg} and $<3, 4, 5> $ for Corollary~\ref{lgmmsgc}.
 \end{remark}
 
 \bigskip
 
Let $\D_b$ be the set of all $b$-digital semigroups
which satisfy the condition stated in Proposition~\ref{dsns}~(iii).
 An inspection of the proof of  Theorem~\ref{ldsgrchar}  immediately yields 
 the following extensions of  the respective results of   
 \cite[Section~2]{rosbrantorr}.
 
 \bigskip
 
\begin{corollary}\label{ldsgrcphi}
The correspondence 
$\theta_b : \LL \to  \D_b$  given
by   $$\theta_b (S):= \set{z \in Z_b  \;:\;   \ell_b (z) \in S}$$
 is a bijective map, and its inverse 
 $\varphi_b : \D_b \to  \LL$  is defined by
$$\varphi_b (D):= L_b (D)  \cup  \set{0}.$$ 
\end{corollary}
 
 \bigskip
 
\begin{corollary}\label{ldsgrcnm}
For every $D \in  \D_b$  the set   $Z_b  \setminus D$ is finite.
\end{corollary}
\begin{proof}
By what we have seen so far we know that $S:=\varphi_b (D)$ is a numerical monoid.
If $b>1$ then analogously as in the proof of \cite[Corollary~8]{rosbrantorr} we show that
 $\set{b^{F(S)}, \to}\subseteq D$.  Now, let  $b<-1$ and $n\in \NN$ be even such that
 $n\ge F(S).$ Then Corollary~\ref{merhrpc}  yields $\set{b^{n}, \to}\subseteq D$. Moreover, 
 $b^{2 n + 1} \in  D$ by  Proposition~\ref{dsns}, hence $(-\infty, \; b^{2 n - 1})  \cap \ZZ  \subseteq D$ by Lemma~\ref{lzzvgl}, and we are done.
 \end{proof}

\bigskip
 
\begin{example}
We have $Z_b \in \D_b$, but  $\NN  \setminus \set{0} \in \D_b$ if and only if  $b>1$.
 \end{example}

\bigskip

Recall that a  $(v, b, r, k)$-configuration is an incidence structure  with $v$ points, $b$ lines, $r$
lines  through each point  and $k$ points on each line.  
Let  $S_{(r, k)}$  be the set of  all  integers $d$ such that there exists a
$(d \cdot \frac{k}{\gcd {(r, k)}}, d \cdot \frac{r}{\gcd {(r, k)}} , r, k)$-configuration.
Bras-Amor\'os and Stokes \cite[Theorem~2]{brasamorosstokes} showed that $S_{(r, k)}$ is a 
numerical monoid provided $r, k \ge 2$. By  \cite[Introduction]{rosbrantorr} 
$S_{(r, r)}$ is an LD-semigroup if $r\ge 2$, and this statement is slightly sharpened now.
 
 \begin{theorem}
 If   $r  \ge 2$  then $S_{(r, r)}$ belongs to $\LL_-$.
\end{theorem} 
\begin{proof}
Let    $S:=S_{(r, r)}$ and $s, t \in S \setminus \set{0,1}$.
By  \cite[Section~2]{stokesbrasamoros} we know that $s+t-1, s+t+1\; \in S$.
Therefore, in view of     Theorem~\ref{ldsgrchar} it suffices to show that 
 $s+t-3\in S$.

If $r=2$ then we infer  $S=<3,4,5>$ from  \cite[Corollary~1]{brasamorosstokes},
and we easily deduce our claim.
Now, let $r>2$ and $m$ be the multiplicity of  $S$, i.e., the least positive integer belonging to $S$. Then we have
$$m\ge r^2 - r +1\ge 3$$
by   \cite[Lemma~1]{stokesbrasamoros}. Since we may assume $s\ge t\ge m$ we
find $s+t-3\in S$ by  \cite[Theorem~9]{stokesbrasamoros}.
  \end{proof}

\section{Generating  $b$-digital semigroups}

This section is devoted to a description of the set $\D_b$ which is very closely related
to the respective result in \cite{rosbrantorr}. Let us start with the analogue of 
\cite[Lemma~16]{rosbrantorr} which can immediately be verified.

\begin{lemma} \label{intsecdsg}
The intersection of  $b$-digital semigroups which belong to $\D_b$ is a
 $b$-digital semigroup in $\D_b$.
    \end{lemma}

In view of this result, given $A\subseteq Z_b$ the set
$$\D_b (A):= \bigcap_{D  \in \D_b, \; A\subseteq D}  D $$
 is the smallest
element of $\D_b$ which contains $A$.

\bigskip

For $A\subseteq \NN \setminus \set{0}$ we let  $\LL_b (A)$ denote the intersection of all
$b$-LD-semigroups which contain $A$. Analogously  as 
\cite[Proposition~17, Corollary~18]{rosbrantorr} 
we write down  the following result based on  Theorem~\ref{ldsgrchar} and  
\cite[Lemma~1]{rosbrantorr}.

\begin{proposition}\label{lbdldsg} 
If $A\subseteq \NN \setminus \set{0}$ is nonempty then  $\LL_b (A)$ is the smallest
$b$-LD-semigroup which contains $A$.
\end{proposition}
 
\bigskip

Now we straightforwardly extend \cite[Proposition~19]{rosbrantorr}.

\begin{proposition}\label{lbdldsg} 
Let $S\in \LL$ and $A\subseteq \NN \setminus \set{0}$ be nonempty. Then  $S$ is the smallest
$b$-LD-semigroup containing $\LL_b (A)$ if and only if   $\theta_b(S)$ is the smallest
element of $\D_b$ which contains~$A$.
\end{proposition}

\bigskip

Let $A$ be a subset of the $b$-digital semigroup $D$. Following \cite[Section~4]{rosbrantorr} we call $A$ a  $\D_b$-system of 
generators of $D$ if  $\D_b (A)=D$; we say that  $A$ is a  minimal $\D_b$-system of 
generators  of $D$ if  no proper subset of $A$  is a  $\D_b$-system of 
generators of $D$. Analogously as  \cite[Theorem~21]{rosbrantorr} we can 
prove the following theorem using  Lemma~\ref{intsecdsg}, Corollary~\ref{ldsgrcphi},
Corollary~\ref{ldsgrcharclds} and  Proposition~\ref{lbdldsg}.
  
\begin{theorem}\label{dbfingen} 
We have
$$\D_b =  \set{\D_b (A)  \;:\; A \text{ finite nonempty subset of }  Z_b}. $$
\end{theorem}

\section{$b$-LD-semigroups containing prescribed integers}

In this section we treat $b$-LD-semigroups which contain a prescribed  set of positive integers.
In particular, we derive an algorithm calculating the smallest element of $\LL_-$
which contains given positive integers. Due to the fact that $E_b$ may contain a positive element we  present  a restricted $b$-adic version of \cite[Proposition~28]{rosbrantorr}.

\begin{proposition}\label{bldsgchar} 
Let $S\ne \NN$ be a numerical monoid and $\msg (S) = \set{n_1, \ldots, n_p}$. Then  the following
statements are equivalent:
\begin{enumerate}
  \item $S$ is a $b$-LD-semigroup.
    \item  If  $e\in E_b$ and $i, j \in  \set{1, \ldots,p}$ then $n_i + n_j +e \in S$.
  \item  If  $e\in E_b$ and $s\in S \setminus  \set{0, n_1, \ldots,n_p}$ then $s +e \in S$.
\end{enumerate}
 \end{proposition}
  \begin{proof}
  (i) $\implies$ (ii): Clear by Theorem~\ref{ldsgrchar}.\\
  (ii) $\implies$ (iii): Let $t\in S$ and $i, j \in  \set{1, \ldots,p}$ such $s=n_i + n_j +t$.
  Then we clearly have $$s+e= (n_i + n_j +e)+t \in S\,.$$
  (iii) $\implies$ (i): Let $s, t\in S \setminus  \set{0,1}$. Then  $s+t  \in   S \setminus  \set{0, n_1, \ldots,n_p}$, hence  $s+t+e \in S$, and we are done  by Theorem~\ref{ldsgrchar}.
 \end{proof}

\bigskip

It does not seem obvious how \cite[Proposition~28~(iv)]{rosbrantorr} can be modified
for a characterization of the semigroups in $\LL_-$. In fact, both numerical monoids
$S:=  \,  <3, 4, 5>$ and $T:=  \,  < 4, 5, 7>$ belong to  $\LL_-$ and satisfy the conditions given
in  Proposition~\ref{bldsgchar} and \cite[Proposition~28~(iv)]{rosbrantorr}. Furthermore,
we have $$3\in S, \quad 3-(P(3)+1)= 1 \notin S, \quad   3-(P(3)-3)=5\in S,$$
but
$$4\in T, \quad  4-(P(4)+1)= 2 \notin T, \quad   4-(P(4)-3)=6\notin T;$$
here we set 
$$P(s):=  \max \set{c_1 +  \cdots + c_p  \;:\;  c_1,   \ldots,  c_p  \in \NN  
  \text{  and  }  s= c_1 n_1 +  \cdots + c_p n_p}$$
where  $s$ is an element of the numerical monoid with minimal
system of generators $\set{n_1,   \ldots,  n_p}.$

On the other hand, for $U:=  \,  < 2, 3> \, \in \LL  \setminus \LL_-$ we have
$$2\in U \quad  \text{and} \quad  2-(P(2)-3),  \, 2-(P(2)-1),  \, 2-(P(2)+1) \, \in U.$$
 
\bigskip
 
Clearly, in view of  Proposition \ref{bldsgchar},  Theorem~\ref{ldsgrchar}  and \cite[Proposition~28~(iv)]{rosbrantorr}
we can immediately formulate the following result.

\begin{proposition}\label{bldsgposchar} 
Let  $S$ be a numerical monoid and $b>1$. Then $S$ is a $b$-LD-semigroup
if and only if $s- \set{0,  \ldots, P(s) -1} \subset S$ for all 
 $s\in S \setminus  \set{0}$.
 \end{proposition}

\bigskip

The algorithm below computes the smallest element of $\LL_-$ containing
a given finite set of integers larger than $1$. After choosing a large heuristic bound
the algorithm closely follows \cite[Algorithm~32]{rosbrantorr}
for the determination of the smallest LD-semigroup containing a set of positive integers, and in view of
our previous results
the  justification of its behavior is analogous to the one in  \cite[Section~5]{rosbrantorr}.
 
\bigskip

	
\begin{algorithm}

\caption{Computation of the smallest element of $\LL_-$ containing given positive integers}

\begin{algorithmic}

\REQUIRE $\text{Non-void finite subset }  A  \subset \NN \setminus \set{0, 1}, \; {\rm bound} \in \mathbb{N}$.

\ENSURE The minimal system of generators of the smallest element of $\LL_-$ containing $A$  \makebox{ or ``{\rm overflow}''}

\STATE $k \gets 0$

\STATE $E \gets \set{-3, -1, 1}$

\REPEAT

\STATE $k \gets k+1$

\STATE $B \gets \msg (A)$

\STATE $A \gets B \cup  \set{x+y+e\;:\; x,y\in B,\; e\in E,\; x+y+e\notin  \,    <B>}$

\UNTIL{$k > {\rm  bound }$ or $B=A$}

\IF{$k > {\rm bound}$}

\RETURN{\rm ``overflow''}

\ELSE

\RETURN{\rm ``Minimal system of generators:'' B}

\ENDIF

\end{algorithmic}

\end{algorithm}


\bigskip

Let us illustrate this algorithm by an easy example.

 \begin{example}
We determine the minimal system of generators
 of  the smallest element $S$ of $\LL_-$ containing $8$.
  Our algorithm requires the following three steps:
 \begin{itemize}
  \item   $B=\set{8}, \;  A=  B \cup \set{13, 15, 17}$ 
  \item   $B=  \set{8, 13, 15, 17}, \;  A=  B \cup \set{18, 20, 22, 27}$
   \item   $B= \set{8, 13, 15, 17, 18, 20, 22, 27},  \;  A=  B $
\end{itemize}
Therefore
$$S=  \,  <{8, 13, 15, 17, 18, 20, 22,  27}>  \,  = \set{0, 8, 13, 15, 16, 17, 18, 20,\to}.$$
It seems worthwile to remark that $S$ is not  an Arf numerical semigroup  (see \cite{rgsggb}),
because $2\cdot 16 - 13 = 19 \notin S$.
 \end{example}

\section{Auxiliary results on the lengths of $b$-adic representations}

The considerations presented in the previous sections are based on the knowledge of 
some facts on the lengths
of $b$-adic representations of integers. These facts are certainly well-known, but are collected here for the sake of completeness. First we recall a fundamental  observation which is tacitly used in this paper.

\begin{proposition}\label{modlen} \cite[Proposition 3.1]{frolai11}
Let $ b < -1 $ and  $z\in \ZZ$.
If $z>0$ then $\ell_b (z)$ is odd, and if $z<0$ then $\ell_b (z)$ is even.
\end{proposition}

\begin{example} \label{exrbbme} 
Let  $b<-1$ and $u\in N_b \setminus  \set{0}$.
  Then  we have $-u = b+v$ with some $v\in N_b$, thus $\ell_b (-u)=2$.
In particular, we have    $-1 = b + (\abs{b}-1)$, hence the base $b$ representation of $\abs{b}$ is 
$$\abs{b}= (-1) \cdot b  = b^2 + (\abs{b}-1)\, b\,,$$
and we have $\ell_b (\abs{b})=3$.
 \end{example}

\bigskip

Using \cite[Lemma 7]{matula} the following bounds for the
length of the $b$-adic representation of an integer~$z$
can immediately be derived:
$$\frac{\log \abs{z} - \log (\abs{b}-1) }{\log \abs{b}} \le  \ell_b (z) \le \frac{\log \abs{z} }{\log \abs{b}}+4 \qquad \qquad (z\in  Z_b)$$

\bigskip

 However, our purposes require  bounds which depend on
the signs of the integers $b$ and $z$.  
Note that the next  result yields
an explicit description of the sets $\Delta_b (n)$.

\bigskip

\begin{proposition} \label{merhrp}
 Let $ b \in \ZZ \setminus \set{-1,0,1}$ and  $a\in \ZZ$.
 \begin{enumerate} 
 \item  If $b > 1$  and $a> 0$ then $ \ell_b (a)=\ell$ if and only if 
 $$ b^{\ell - 1}\le a \le b^{\ell }-1. $$
In this case we have 
$$\frac{\log a}{\log b} < \ell \le \frac{\log a}{\log b}+1 \,.$$
  \item  If $b <-1$ and $a>0$ then  $ \ell_b (a)=\ell$ if and only if 
 $$ \frac{b(b^{\ell - 2}-1)}{1-b}\le a \le \frac{b^{\ell +1}-1}{1-b}\,.$$
In this case we have
$$\frac{\log \bigl( (\abs{b} +1)a +1\bigr)}{\log \abs{b}} -1 \le  \ell \le \frac{\log \bigl( (1+ 1/ \abs{b} )a -1\bigr)}{\log \abs{b}}+2 \,.$$
  \item  If $b <-1$ and $a<0$ then $ \ell_b (a)=\ell$ if and only if 
 $$ \frac{b(b^{\ell }-1)}{1-b}\le a \le \frac{b^{\ell -1}-1}{1-b}\,.$$
In this case we have
$$\frac{\log \bigl( (1+ 1/\abs{b} ) \abs{a} +1\bigr)}{\log \abs{b}} \le  \ell \le \frac{\log \bigl( (1+  \abs{b} ) \abs{a} -1\bigr)}{\log \abs{b}}+1 \,.$$
 \end{enumerate}
 \end{proposition}
\begin{proof}
(i) This is well-known and easy to check.\\
(ii) We observe
$$a\le (\abs{b}-1) \sum_{i=0}^{(\ell -1)/2}  b^{2i}= -(b+1)  \frac{b^{2((\ell -1)/2+1)}-1}{b^2-1}= - \frac{b^{\ell+1}-1}{b-1}$$
and
$$a\ge b^{2 \cdot (\ell -1)/2} + (\abs{b}-1) \sum_{i=1}^{(\ell -1)/2}  b^{2i-1}$$
$$=b^{\ell -1} - \frac{b+1}{b}\Bigl(\frac{b^{2((\ell -1)/2+1)}-1}{b^2-1}-1\Bigr)=  \frac{b(1-b^{\ell-2})}{b-1}\,,$$
from which the estimates for $\ell$ are derived straightforwardly.\\
(iii) Noting 
$$ (\abs{b}-1)\bigl(b^{2\cdot \ell /2-1}+  \sum_{i=1}^{\ell /2-1}  b^{2i-1}\bigr) \le a \le b^{\ell -1}+   (\abs{b}-1)\sum_{i=0}^{\ell /2-1}  b^{2i}$$
we complete the proof  as above.
 \end{proof}

\begin{corollary} \label{merhrpc}
 Let $ b < -1$ and  $a, n \in \NN$. If $n$ is even and $a\ge b^n$ then we have $ \ell_b (a)>n$.
 \end{corollary}
\begin{proof}
Assume the contrary. Then Proposition \ref{modlen} yields $ \ell_b (a) \le n-1$, hence
$n\ge 2$ and we infer
the impossible inequality $b^n \le (b^n-1)/ (1-b)$
from the Proposition.
 \end{proof}
 
\bigskip

Now we compare the sizes of integers to the lengths of their $b$-adic
representation.

\begin{lemma} \label{lzzvgl}
 Let 
 $a, c \in \ZZ$.
 \begin{enumerate} 
 \item   If $0 \le a < c$ then $  \ell_b (a)\le \ell_b (c).$
\item   If $ a,  c\ge 0 \text{ and }  \ell_b (a) <  \ell_b (c)$ then $  a < c.$
 \item  If $b <-1$ and $a>0$ then we have $ \ell_b (-a)=\ell_b (a)+1$.
\item  Let $b <-1$. 
 \begin{enumerate} 
\item   $ a < c<0 \implies  \ell_b (a)\ge \ell_b (c).$
\item   $ a,  c\le 0 \text{ and }  \ell_b (a) >  \ell_b (c)  \implies  a < c.$
  \end{enumerate}
 \end{enumerate}
 \end{lemma}
\begin{proof}
(i) This is well-known and easy to check.\\
(ii) -- (iv) This is straightforwardly derived from Proposition~\ref{merhrp}.
 \end{proof}

\begin{lemma} \label{znlm}
 Let $ b < -1$ and  $n, m$ be even positive integers such that $n\le m$.
 If  $$\frac{b (b^{m}-1)}{1-b}   \le z   \le \frac{b^{n-1}-1}{1-b}$$
 then we have  $$ n\le \ell_b (z)\le m\,.$$
 \end{lemma}
\begin{proof}
Let  $y \in \ZZ$ such that  $ \ell_b (y) = n$ and assume 
$ n >   \ell_b (z)\,.$ Then  Proposition~\ref{merhrp} and Lemma~\ref{lzzvgl} yield
$$\frac{b (b^{n}-1)}{1-b}   \le y<z   \le \frac{b^{n-1}-1}{1-b}$$
and then $ n =  \ell_b (z)\,:$ Contradiction.

The second inequality is proved analogously.
 \end{proof}

 \bigskip
 
Further, we need the length of the $b$-adic representation of the  product of two elements. 
 
\begin{lemma} \label{lacn}
  \begin{enumerate}
  \item	Let  $ b > 1$ and  $a, c \in \NN\setminus \set{0}$. Then we have 
  \begin{equation}\label{lacneq}
  \ell_b (a c) =   \ell_b (a )+  \ell_b ( c)+ e
\end{equation}
 for some  $e \in \set{-1, 0}. $
   \item	Let $ b < -1$ and  $a, c \in \ZZ \setminus  N_b$. Then 
   there exists some $e \in \set{-3,-1, 1}$  such that  \eqref{lacneq} holds. 
\item If $n, m \ge 2$ and $e \in E_b$  
then there exist
    $a, c \in Z_b$ such that $\ell_b (a )=n, \;  \ell_b ( c)=m$ and   \eqref{lacneq} holds.
         \end{enumerate}
 \end{lemma}
\begin{proof}
(i) For $b=2$ this is immediately checked using Proposition~\ref{merhrp}, and for $b>2$ the proof of 
\cite[Lemma~3]{rosbrantorr} can easily be extended.\\
(ii) Set $n:=   \ell_b (a )$ and $m:= \ell_b ( c)$. Certainly it suffices to consider the subsequent cases.

\medskip

Case 1 \hspace{3cm}  $a>0$

\medskip

Then $n$ is odd and we infer
 $$ \frac{b(b^{n - 2}-1)}{1-b}\le a \le \frac{b^{n +1}-1}{1-b}$$
from  Proposition~\ref{merhrp}. 

\medskip

Case 1.1  \hspace{27mm}  $c>0$

\medskip

Then $m$ is odd and as above we have
 $$ \frac{b(b^{m - 2}-1)}{1-b}\le c \le \frac{b^{m +1}-1}{1-b}\,.$$
Now we easily verify
 $$ \frac{b(b^{n+m - 5}-1)}{1-b}\le \frac{b^2(b^{n+m - 4}-b^{n - 2}-b^{m - 2}+1)}{(1-b)^2}\le a c $$
$$\le \frac{b^{n+m +2}-b^{n +1}-b^{m +1}+1}{(1-b)^2}\le \frac{b^{n+m +2}-1}{1-b}\,.$$
Then Proposition~\ref{merhrp} yields 
$$ n+m-3\le \ell_b (a c)  \le n+m+1\,,$$
and our assertion follows from Proposition \ref{modlen}.

\medskip

Case 1.2  \hspace{27mm}  $c<0$

\medskip

As above we  verify
 $$ \frac{b(b^{n+m +1}-1)}{1-b}\le \frac{b(b^{n+m +1}-b^{n+1}-b^{m }+1)}{(1-b)^2}\le a c $$
$$\le \frac{b^{n+m -3}-b^{n -2}-b^{m -1}+1}{(1-b)^2}\le \frac{b^{n+m -4}-1}{1-b}$$
 keeping in mind that $m$ is even, and then we conclude using 
 Lemma~\ref{znlm}.

\medskip

Case 2  \hspace{3cm}  $a<0$

\medskip

We may suppose $c<0$ and proceed as in Case 1.1\,.

(iii) 
The case  $ b >1$ is well-known. Now, let  $ b <-1$.
   For the positive integers  $$a= \frac{b^{2n}-1}{1-b}   \quad  \text{ and  } \quad  c= \frac{b^{2m}-1}{1-b}$$
  we have $$ \ell_b (a c) =   \ell_b (a )+  \ell_b ( c) + 1\,.$$
 Similarly, for the negative integers  $$a= \frac{b^{2n  -1}-1}{1-b}   \quad  \text{ and  } \quad  c= \frac{b (b^{2m}-1)}{1-b}$$
  we verify  $$ \ell_b (a c) =   \ell_b (a )+  \ell_b ( c) - 1\,,$$
 and for   $$a= \frac{b (b^{2n  -3}-1)}{1-b}   \quad  \text{ and  } \quad  c= \frac{b (b^{2m-1}-1)}{1-b}$$
  we see $$ \ell_b (a c) =   \ell_b (a )+  \ell_b ( c) - 3\,.$$
    \end{proof}

 \bigskip

We close this section by an easy application of Proposition~\ref{merhrp}
the details of which we leave to the reader  (cf. the special case $b=10$ in \cite[proof of Corollary~9]{rosbrantorr}).

\begin{proposition} \label{cdbn}
 For  
  $n\in \NN\setminus \set{0}$
 we  have 
$$ \Card \, (\Delta_b (n)) =
\begin{cases}  (b-1)\,  b^{n-1} & (b > 1),\\
-   (b+1)\, b^{n-1}  & (b < - 1, \; n  \text{ odd}),\\
   (b+1)\,b^{n-1} & (b < - 1, \; n  \text{ even}).\\
\end{cases}
$$
 \end{proposition}

\bigskip

{\bf Acknowledgement}.
The author  is indebted to Denise  Torr{\~a}o for bringing the work \cite{rosbrantorr} to his knowledge and to anonymous referees for very carefully 
reading the first version of this paper.

\bigskip

\def\cprime{$'$}
\bibliographystyle{siam}

\end{document}


\vspace{-10mm}

\centerline{{\bf Corrigendum}}

\medskip

\centerline{ to the article}

\bigskip

\centerline{{ Horst {Brunotte}}}

\centerline{{\em {Digital semigroups}}}

   \centerline{{RAIRO-Theor. Inf. Appl.}  50  (2016) 67 -- 79.}

\bigskip

\date{\today}

\bigskip

p. 70: Replace  the first part of the proof of Theorem 3.2 by:

\medskip

(i) $\implies$ (ii) 
Let  $D$ be a $b$-digital semigroup such that  $S= L_b (D) \cup  \set{0}$. Then we clearly have
 $S\ne  \set{0}$.  Let $s, t  \in S\setminus \set{0, 1}$ and $e \in E_b \,.$ W.l.o.g. we may assume
$2\le s \le t$. 

First, consider  $s\ge 3$. By  Lemma~6.7 there exist $a, c \in Z_b$ such that $s= \ell_b (a), \, t= \ell_b (c)$
 and $\ell_b (a c) = s + t +  e \,.$
 By the properties of  $D$ we know that  $a, c \in D$, thus
$$ s+t+ e = \ell_b (a c)  \in  L_b (D) \subset S\,.$$

Second, let  $s=2$. For $b\ne -2$ we conclude as above. For $b =  -2$ we observe
 $\ell_{-2} (-1)=2$, hence $-1\in D$, thus $1=(-1)^2\in D$ and $ 1 = \ell_{-2}  (1)  \in  S$
yielding $S=\NN$, and our claim drops out.

\bigskip

p. 75 -- 76: Replace  Proposition 6.3 (ii) by:

\medskip

If $b <-1$ and $a>0$ then  $ \ell_b (a)=\ell$ if and only if 
 $$ \frac{b^{\ell}-b^2+1}{b(1-b)}\le a \le \frac{b^{\ell +1}-1}{1-b}\,.$$
In this case we have
$$\frac{\log \bigl( (\abs{b} +1)a +1\bigr)}{\log \abs{b}} -1 \le  \ell \le \frac{\log \bigl( \abs{b} (\abs{b}+1)a- b^2+1\bigr)}{\log \abs{b}} \,.$$
   
\bigskip

p. 76: Replace  the last equality in the proof of Proposition 6.3 (ii) by 
 $$ =\frac{b^{\ell}-b^2+1}{b(1-b)}$$

\bigskip

p. 76: Replace Lemma 6.5 (iii) by:

\medskip

 If $b <-1$ and $a>0$ then we have $ \ell_b (-a) \in \set{\ell_b (a)-1, \ell_b (a)+1}$.

\bigskip

p. 77:  Delete Lemma 6.6  since it is not needed.

\bigskip

p. 77: Replace  Lemma 6.7 (iii) by:

\medskip

If $n, m \ge 2$ ($n, m \ge 3$ in case $b=-2$) and $e \in E_b$  
then there exist
    $a, c \in Z_b$ such that $\ell_b (a )=n, \;  \ell_b ( c)=m$ and  
(6.1) holds.

\bigskip

p. 77 -- 78: Replace the proof of Lemma 6.7 (ii) by:

\medskip

Set $n:=   \ell_b (a )$ and $m:= \ell_b ( c)$. We consider the subsequent cases.

\medskip

Case 1 \hspace{3cm}  $a>0$

\medskip

Then $n$ is odd and we infer
 $$  \frac{b^{n}-b^2+1}{b(1-b)}   \le a \le \frac{b^{n +1}-1}{1-b}$$
from  Proposition~6.3.   

\medskip

Case 1.1  \hspace{27mm}  $c>0$

\medskip

Then $m$ is odd and we have
 $$ \frac{b^{m}-b^2+1}{b(1-b)}   \le c \le \frac{b^{m +1}-1}{1-b}\,.$$
In view of 
$$ \frac{b^{n+m - 3}-b^2+1}{b(1-b)} \le  a c \le  \frac{b^{n+m +2}-1}{1-b}$$ 
Proposition~6.3 
yields 
$$ n+m-3\le \ell_b (a c)  \le n+m+1      \qquad \qquad  \qquad \qquad   (6.2)$$
and our assertion follows from Proposition 6.1. 

\medskip

Case 1.2  \hspace{27mm}  $c<0$

\medskip

Then $m$ is even, we have
$$ \frac{b(b^{m }-1)}{1-b}\le c \le \frac{b^{m -1}-1}{1-b}$$
and
$$ \frac{b(b^{n+1 }-1)(b^{m }-1)}{(1-b)^2}\le a c \le \frac{(b^{n} -b^2+1)(b^{m -1}-1)}{b(1-b)^2}\,.$$
We convince ourselves that (6.2) holds, 
and an application of Proposition 6.1 terminates the proof of this case.

\medskip

Case 2  \hspace{3cm}  $a<0$

In view  of what we have seen so far, we only need to treat the case  $c<0$ and here we proceed analogously as in Case 1.1.

\bigskip

p. 78: Replace the proof of Lemma 6.7 (iii) by:

\medskip

The case  $ b >1$ is well-known. Now, let  $ b <-1$ and exploit Propositions 6.1 and 6.3.

First, let $n$ be odd. If  $m$ is odd we choose
$$ a =  \frac{b^{n +1}-1}{1-b}  \qquad \text{ and } \qquad c =  \frac{b^{m +1}-1}{1-b} \qquad (e=1),$$
$$ a =  \frac{b^{n +1}-1}{1-b}  \qquad \text{ and } \qquad c =  \frac{b^{m }-b^2 +1}{b(1-b)} \qquad (e=-1),$$
and 
$$ a =  \frac{b^{n}-b^2 +1}{b(1-b)} \qquad \text{ and } \qquad c =  \frac{b^{m }-b^2 +1}{b(1-b)} \qquad (e=-3),$$
and if $m$ is even we take
$$ a =     \frac{b^{n +1}-1}{1-b}  \qquad \text{ and } \qquad c =  \frac{b (b^{m}-1)
}{1-b} \qquad (e=1),$$
$$ a =    \frac{b^{n}-b^2 +1}{b(1-b)}  \qquad \text{ and } \qquad c =  \frac{b (b^{m } - 1)}{1-b} \qquad (e=-1),$$
and 
$$ a =  \frac{b^{n}-b^2 +1}{b(1-b)} \qquad \text{ and } \qquad c =  \frac{b^{m-1 }- 1}{1-b} \qquad (e=-3).$$

\medskip

Second, for  $n, m$ even we consider
$$ a =  \frac{b (b^{n } - 1)}{1-b} \qquad \text{ and } \qquad c = \frac{b (b^{m} - 1)}{1-b} \qquad (e=1),$$
$$ a =  \frac{b (b^{n}- 1)}{1-b} \qquad \text{ and } \qquad c = \frac{b^{m -1} - 1}{1-b} \qquad (e=-1),$$
and
$$ a =  \frac{b^{n -1} - 1}{1-b} \qquad \text{ and } \qquad c = \frac{b^{m -1} - 1}{1-b} \qquad (e=-3).$$

\bigskip